\documentclass[12pt]{article}%
\usepackage{graphicx}
\usepackage[intlimits]{amsmath}
\usepackage{latexsym}
\usepackage{amsfonts}
\usepackage{amssymb}%
 \makeatletter
\newfont{\footsc}{cmcsc10 at 8truept}
\newfont{\footbf}{cmbx10 at 8truept}
\newfont{\footrm}{cmr10 at 10truept}
\newtheorem{theorem}{Theorem}

\newtheorem{lemma}[theorem]{Lemma}

\newenvironment{proof}[1][Proof]{\noindent{\textbf {#1}  }}  {\hfill$\Box$\bigskip}

\begin{document}

\title{Graphs with many $r$-cliques have large complete $r$-partite subgraphs}
\author{Vladimir Nikiforov\\{\small Department of Mathematical Sciences, University of Memphis, Memphis TN
38152}}
\maketitle

\begin{abstract}
Let $r\geq2$ and $c>0.$ Every graph on $n$ vertices with at least $cn^{r}$
cliques on $r$ vertices contains a complete $r$-partite subgraph with $r-1$
parts of size $\left\lfloor c^{r}\log n\right\rfloor $ and one part of size
greater than $n^{1-c^{r-1}}.$ This result implies the
Erd\H{o}s-Stone-Bollob\'{a}s theorem, the essential quantitative form of the
Erd\H{o}s-Stone theorem.\medskip


\end{abstract}

\subsection*{Main results}

This note is part of an ongoing project aiming to renovate some classical
results in extremal graph theory, see, e.g., \cite{BoNi04},
\cite{Nik07,Nik07c}.

Recall a result of Erd\H{o}s and Stone \cite{ErSt46}: every graph with $n$
vertices and $\left\lceil cn^{2}\right\rceil $ edges contains a complete
bipartite subgraph with each part of size $\left\lfloor a\log n\right\rfloor
$, where $c,a>0$ are independent of $n.$

One generalization of this result stems from the fundamental theorem of
Erd\H{o}s and Stone \cite{ErSt46}: every graph with $n$ vertices and
$\left\lceil \left(  1-1/r+c\right)  n^{2}/2\right\rceil $ edges contains a
complete $\left(  r+1\right)  $-partite graph with each part of size $g\left(
r,c,n\right)  ,$ where, for $r$ and $c$ fixed, $g\left(  r,c,n\right)  $ tends
to infinity with $n$. In \cite{BES76} Bollob\'{a}s and Erd\H{o}s found that
$g\left(  r,c,n\right)  =\Theta\left(  \log n\right)  ,$ and in \cite{BoEr73},
\cite{BoKo94}, \cite{ChSz81}, and \cite{Ish02} the function $g\left(
r,c,n\right)  $ was determined with great precision.

In this note we deduce the Erd\H{o}s-Stone-Bollob\'{a}s result from weaker
premises. Also, in our setup $c$ may be a function of $n,$ say $c=1/\ln\ln n$
and the bounds on $g\left(  r,c,n\right)  $ remain meaningful; such results
are beyond the scope of the papers mentioned above.

Our notation follows \cite{Bol98}. Thus, $K_{r}\left(  s_{1},\ldots
,s_{r}\right)  $ denotes a complete $r$-partite graph with parts of size
$s_{1},\ldots,s_{r},$ and an $r$-clique is a complete subgraph on $r$
vertices. We write $E\left(  G\right)  $ for the edge set of a graph $G$ and
$K_{r}\left(  G\right)  $ for the set of its $r$-cliques; we let $e\left(
G\right)  =\left\vert E\left(  G\right)  \right\vert $ and $k_{r}\left(
G\right)  =\left\vert K_{r}\left(  G\right)  \right\vert $.

Given $M\subset K_{r}\left(  G\right)  ,$ we write $K_{s}\left(  M\right)  $
for the set of $s$-cliques contained in members of $M.$ Given a subgraph
$H\subset G$ such that $H=K_{r}\left(  s_{1},\ldots,s_{r}\right)  $, we say
that $M$ \emph{covers} $H$ if $E\left(  H\right)  \subset K_{2}\left(
M\right)  $ and $H$ contains $\min\left\{  s_{1},\ldots,s_{r}\right\}  $
disjoint members of $M.$

Here is our main result.

\begin{theorem}
\label{genZ}Let $r\geq2,$ $c^{r}\ln n\geq1,$ and $G$ be a graph with $n$
vertices. Every set of at least $cn^{r}$ $r$-cliques of $G$ covers a
$K_{r}\left(  s,\ldots s,t\right)  $ with $s=\left\lfloor c^{r}\ln
n\right\rfloor $ and $t>n^{1-c^{r-1}}.$
\end{theorem}

\subsubsection*{Remarks}

\begin{itemize}
\item[-] Using random graphs, it is easy to see that most graphs on $n$
vertices contain no complete bipartite subgraphs with both parts larger than
$C\log n,$ for some $C>0,$ independent of $n.$ Hence, Theorem \ref{genZ} is
essentially best possible.

\item[-] Theorem \ref{genZ} implies the Erd\H{o}s-Stone-Bollob\'{a}s theorem.
Indeed, in \cite{KhNi78} (see also \cite{Lov79}, Problem 11.8) it is proved
that if $k_{s}\left(  G\right)  >0,$ then
\[
\frac{\left(  s+1\right)  k_{s+1}\left(  G\right)  }{sk_{s}\left(  G\right)
}-\frac{n}{s}\geq\frac{sk_{s}\left(  G\right)  }{\left(  s-1\right)
k_{s-1}\left(  G\right)  }-\frac{n}{s-1}.
\]
Thus, if $e\left(  G\right)  \geq\left(  1-1/r+c\right)  n^{2}/2,$ we find
that $k_{r+1}\left(  G\right)  >\left(  c/r^{r}\right)  n^{r+1},$ and so $G$
contains a $K_{r+1}\left(  s,\ldots s,t\right)  $ with
\[
s=\left\lfloor \left(  c/r^{r}\right)  ^{r+1}\ln n\right\rfloor \text{,
\ \ }t>n^{1-\left(  c/r^{r}\right)  ^{r}}.
\]
This is slightly stronger than the Erd\H{o}s-Stone-Bollob\'{a}s result and is
comparable with those established in \cite{BoKo94}.

\item[-] The exponent $1-c^{r-1}$ in Theorem \ref{genZ} is not the best one,
but is simple.\bigskip
\end{itemize}

We deduce Theorem \ref{genZ} from a routine lemma, proved here for convenience.

\begin{lemma}
\label{le1} Let $F$ be a bipartite graph with parts $A$ and $B.$ Let
$\left\vert A\right\vert =m,$ $\left\vert B\right\vert =n,$ $r\geq2,$ $\left(
\ln n\right)  ^{-1/r}\leq c<1/2,$ and $s=\left\lfloor c^{r}\ln n\right\rfloor
.$ If $s\leq\left(  c/2\right)  m+1$ and $e\left(  F\right)  \geq cmn,$ then
$F$ contains a $K_{2}\left(  s,t\right)  $ with parts $S\subset A$ and
$T\subset B$ such that $\left\vert S\right\vert =s$ and $\left\vert
T\right\vert =t>n^{1-c^{r-1}}$.
\end{lemma}

\begin{proof}
Let
\[
t=\max\left\{  x:\text{there exists }K_{2}\left(  s,x\right)  \subset F\text{
with part of size }s\text{ in }A\right\}  .
\]

For any $X\subset A,$ write $d\left(  X\right)  $ for the number of vertices
joined to all vertices of $X.$ By definition, $d\left(  X\right)  \leq t$ for
each $X\subset A$ with $\left\vert X\right\vert =s;$ hence,%
\begin{equation}
t\binom{m}{s}\geq%
%TCIMACRO{\tsum \limits_{X\subset A,\left\vert X\right\vert =s}}%
%BeginExpansion
{\textstyle\sum\limits_{X\subset A,\left\vert X\right\vert =s}}
%EndExpansion
d\left(  X\right)  =%
%TCIMACRO{\tsum \limits_{u\in B}}%
%BeginExpansion
{\textstyle\sum\limits_{u\in B}}
%EndExpansion
\binom{d\left(  u\right)  }{s}. \label{in1}%
\end{equation}
Set
\[
f\left(  x\right)  =\left\{
\begin{array}
[c]{cc}%
\binom{x}{s} & \text{if }x\geq s-1\\
0 & \text{if }x<s-1,
\end{array}
\right.
\]
and note that $f\left(  x\right)  $ is a convex function. Therefore,%
\[%
%TCIMACRO{\tsum \limits_{u\in B}}%
%BeginExpansion
{\textstyle\sum\limits_{u\in B}}
%EndExpansion
\binom{d\left(  u\right)  }{s}=%
%TCIMACRO{\tsum \limits_{u\in B}}%
%BeginExpansion
{\textstyle\sum\limits_{u\in B}}
%EndExpansion
f\left(  d\left(  u\right)  \right)  \geq nf\left(  \frac{1}{n}%
%TCIMACRO{\tsum \limits_{u\in B}}%
%BeginExpansion
{\textstyle\sum\limits_{u\in B}}
%EndExpansion
d\left(  u\right)  \right)  =n\binom{e\left(  F\right)  /n}{s}\geq n\binom
{cm}{s}.
\]
Combining this inequality with (\ref{in1}), and after rearranging, we find
that%
\begin{align*}
t  &  \geq n\frac{cm\left(  cm-1\right)  \cdots\left(  cm-s+1\right)
}{m\left(  m-1\right)  \cdots\left(  m-s+1\right)  }>n\left(  \frac{cm-s+1}%
{m}\right)  ^{s}\geq n\left(  \frac{c}{2}\right)  ^{s}\\
&  \geq n\left(  e^{\ln\left(  c/2\right)  }\right)  ^{c^{r}\ln n}%
=n^{1+c^{r}\ln\left(  c/2\right)  }=n^{1+2c^{r-1}\left(  c/2\right)
\ln\left(  c/2\right)  }.
\end{align*}
Since $c/2<1/4<1/e$ and $x\ln x$ is decreasing for $0<x<1/e,$ we see that%
\[
t>n^{1+2c^{r-1}\left(  c/2\right)  \ln c/2}>n^{1+c^{r-1}\left(  1/2\right)
\ln1/4}>n^{1-c^{r-1}},
\]
completing the proof.\bigskip
\end{proof}

\begin{proof}
[\textbf{Proof of Theorem \ref{genZ}}]Let $M\subset K_{r}\left(  G\right)  $
satisfy $\left\vert M\right\vert \geq cn^{r}.$ To prove that $M$ covers a
$K_{r}\left(  s,\ldots s,t\right)  $ with $s=\left\lfloor c^{r}\ln
n\right\rfloor $ and $t>n^{1-c^{r-1}},$ we use induction on $r.$

Suppose $r=2$ and let $A$ and $B$ be two disjoint copies of $V\left(
G\right)  .$ Define a bipartite graph $F$ with parts $A$ and $B,$ joining
$u\in A$ to $v\in B$ if $uv\in M.$ Set $s=\left\lfloor c^{2}\ln n\right\rfloor
$ and note that $s\leq\left(  c/2\right)  n+1.$ Since $e\left(  F\right)
=2\left\vert M\right\vert >cn^{2},$ Lemma \ref{le1} implies that $F$ contains
a $K\left(  s,t\right)  $ with $t>n^{1-c}.$ Hence $M$ covers a $K_{2}\left(
s,t\right)  ,$ proving the assertion for $r=2.$ Assume the assertion true for
$2\leq r^{\prime}<r.$

For a set $N\subset K_{r}\left(  G\right)  $ and a clique $R\in K_{r-1}\left(
N\right)  $ let $d_{N}\left(  R\right)  $ be the number of members of $N$
containing $R$. We first show that there exists $L\subset M$ with $\left\vert
L\right\vert >\left(  c/2\right)  n^{r}$ such that $d_{L}\left(  R\right)
>cn$ for all $R\in K_{r-1}\left(  L\right)  .$ Indeed, set $L=M$ and apply the
following procedure.\medskip

\textbf{While }\emph{there exists an }$R\in K_{r-1}\left(  L\right)  $\emph{
with }$d_{L}\left(  R\right)  \leq cn$ \textbf{do}

\qquad\emph{Remove from }$L$\emph{ all }$r$\emph{-cliques containing }%
$R.$\medskip

When this procedure stops, we have $d_{L}\left(  R\right)  >cn$ for all $R\in
K_{r-1}\left(  L\right)  $ and also
\[
\left\vert M\right\vert -\left\vert L\right\vert \leq cn\left\vert
K_{r-1}\left(  M\right)  \right\vert \leq cn\binom{n}{r-1}<\frac{c}{\left(
r-1\right)  !}n^{r}\leq\frac{c}{2}n^{r},
\]
implying that $\left\vert L\right\vert >\left(  c/2\right)  n^{r},$ as claimed.

Since $K_{r-1}\left(  L\right)  \subset K_{r-1}\left(  G\right)  $ and
\[
\left\vert K_{r-1}\left(  L\right)  \right\vert \geq r\left\vert L\right\vert
/n>r\left(  c/2\right)  n^{r}/n=\left(  rc/2\right)  n^{r-1},
\]
the induction assumption implies that $K_{r-1}\left(  L\right)  $ covers an
$H=K_{r-1}\left(  m,\ldots,m\right)  $ with $m=\left\lfloor \left(
rc/2\right)  ^{r-1}\ln n\right\rfloor .$

For a subgraph $X\subset G$ and a vertex $v\in V\left(  G\right)  ,$ write
$X+v$ for the subgraph of $G$ induced by $V\left(  X\right)  \cup\left\{
v\right\}  $. Let $A$ be a set of $m$ disjoint $\left(  r-1\right)  $-cliques
in $H$ that are members of $K_{r-1}\left(  L\right)  .$ Define a bipartite
graph $F$ with parts $A$ and $B=V\left(  G\right)  ,$ joining $R\in A$ to
$v\in B$ if $R+v\in L.$ Set $s=\left\lfloor c^{r}\ln n\right\rfloor .$ Since
$d_{L}\left(  R\right)  >cn$ for all $R\in K_{r-1}\left(  L\right)  ,$ we have
$e\left(  F\right)  >cmn.$ Also, we find that%
\[
s\leq c^{r}\ln n<\left(  c/2\right)  \left(  rc/2\right)  ^{r-1}\ln n<\left(
c/2\right)  \left\lfloor \left(  rc/2\right)  ^{r-1}\ln n\right\rfloor
+1=\left(  c/2\right)  m+1.
\]

Lemma \ref{le1} implies that there exists $K_{2}\left(  s,t\right)  \subset F$
with parts $S\subset A$ and $T\subset B$ such that $\left\vert S\right\vert
=s$ and $\left\vert T\right\vert =t>n^{1-c^{r-1}}.$ Let $H^{\ast}$ be the
subgraph of $H$ induced by the union of the members of $S;$ clearly, $H^{\ast
}=K_{r-1}\left(  s,\ldots,s\right)  $. Since $R+v\in L$ for all $v\in T$ and
$R\in K_{r-1}\left(  H^{\ast}\right)  ,$ we see that $L$ covers a
$K_{r}\left(  s,\ldots,s,t\right)  $, completing the induction step and the proof.
\end{proof}

\subsubsection*{Concluding remarks}

Finally, a word about the project mentioned in the introduction: in this
project we aim to give wide-range results that can be used further, adding
more integrity to extremal graph theory.\bigskip

\textbf{Acknowledgement} Thanks to B\'{e}la Bollob\'{a}s, Sasha Kostochka, and
Cecil Rousseau for helpful discussions.

\end{document}